\documentclass[11pt]{article}

\input{amssym.def}
\input{amssym}
\usepackage{eufrak}

\newtheorem{theorem}{Theorem}

\newtheorem{definition}[theorem]{Definition}

\newtheorem{lemma}[theorem]{Lemma}
\newtheorem{proposition}[theorem]{Proposition}
\newtheorem{remark}[theorem]{Remark}

\def\beq{\begin{equation}}
\def\be{\begin{equation}}
\def\eeq{\end{equation}}
\def\ee{\end{equation}}
\def\K{{\Bbb K}}
\def\paq{{\pa_q}}
\def\pah{{\pa_\h}}

\def\C{{\Bbb C}}
\def\ga{{\gamma}}
\def\R{{\Bbb R}}

\def\D{{\frak D}}
\def\B{{\cal B}}

\def\N{{\cal N}}
\def\A{{\cal A}}

\def\D{{\cal D}}

\def\Wq{{\cal W}_q}

\def\gg{\mathfrak{g}}
\def\pa{\partial}

\def\W{{\cal W}}
\def\K{\Bbb K}
\def\Z{\Bbb Z}

\def\ot{\otimes}
\def\h{\hbar}

\def\hh{\hbar}

\def\gh{\mathfrak{g}_h}

\def\vv{V^{\ot 2}}

\def\tri{\triangleright}

\def\Ren{R_{\End}}
\def\span{{\rm span}}
\def\Sym{{\rm Sym}}

\def\End{{\rm End}}

\def\De{{\Delta}}
\def\qq{{q^{-1}}}
\def\de{{\delta}}

\def\al{\alpha}

\def\Om{{\Omega}}

\def\A{{\cal A}}
\def\B{{\cal B}}

\def\Om{\Omega}
\def\om{\omega}

\begin{document}

\makeatletter
\renewcommand{\theequation}{{\thesection}.{\arabic{equation}}}
\@addtoreset{equation}{section}
\makeatother

\title{Derivatives in noncommutative calculus and deformation property of
quantum algebras}
\author{
\rule{0pt}{7mm} Dimitri
Gurevich\thanks{gurevich@univ-valenciennes.fr}\\
{\small\it LAMAV, Universit\'e de Valenciennes,
59313 Valenciennes, France}\\
\rule{0pt}{7mm} Pavel Saponov\thanks{Pavel.Saponov@ihep.ru}\\
{\small\it
National Research University Higher School of Economics,}\\
{\small\it International Laboratory of Representation Theory and Mathematical Physics}\\
{\small\it 20 Myasnitskaya Ulitsa, Moscow 101000, Russia}\\
{\footnotesize\it \&}\\
{\small\it IHEP, Division of Theoretical Physics, 142281
Protvino, Russia} }

\maketitle

\begin{abstract}
The aim of the paper is twofold. First, we introduce analogs of (partial) derivatives on certain noncommutative algebras,
including some enveloping algebras and their "braided counterparts" --- the so-called modified Reflection Equation algebras.
With the use of the mentioned derivatives we construct an analog of the de Rham complex on these algebras. Second, we discuss
deformation property of some quantum algebras and show that contrary to a commonly held view, in the so-called q-Witt
algebra there is no analog of the PBW property. In this connection, we discuss different forms of the Jacobi condition related
to quadratic-linear algebras.
\end{abstract}

{\bf AMS Mathematics Subject Classification, 2010:} 46L65, 46L87, 81T75

{\bf Key words:}   Jackson derivative,  q-Witt algebra, Jackson $sl(2)$ algebra,
Jacobi condition,  deformation property, (modified) Reflection Equation algebra

\section{Introduction}

In our recent publications \cite{GPS2, GS3, GS4}  we introduced the notion of partial derivatives on some noncommutative (NC)
algebras, in particular, on the enveloping algebras of the Lie algebras $gl(m)$ and their super- and {\em braided} (see below) analogs. 
These partial derivatives differ from their classical counterparts by the form of the Leibniz rule.

In this connection a natural question arises: given a NC algebra $A$, what operators acting on this algebra can be considered as an
appropriate analogs of partial derivatives? This question is pertinent  if $A$ is a deformation (quantization) of the symmetric algebra
$\Sym(V)$ of a vector space $V$ or its super or braided analog. In this paper we give an answer to this question for the enveloping
algebras of some Lie algebras.

Note that the answer depends on a given Lie algebra. Nevertheless, once such partial derivatives are introduced, we are able to define
an analog of the de Rham complex on the corresponding enveloping algebra. Compared with all known approaches to the
problem of defining such a complex, our method leads to objects possessing good deformation property\footnote{For finitely generated
quadratic-linear(-constant) algebras (in particular, enveloping ones) we deal with, this property means that an analog  of the PBW theorem
is valid for them and homogeneous components of the corresponding  quadratic  algebra have stable dimensions (at least for a generic
value of the  parameter). Note that nowadays the term "PBW property" is often used as a synonym of our "good deformation property".
We prefer to reserve  this term for a deformation of  quadratic algebras by linear(-constant) terms.}. However, the  terms of our complex
are endowed with {\em one-sided} $A$-module structures, whereas  using the classical Leibniz rule for the de Rham operator requires a
two-sided $A$-module structure. Also, we define the notion of the Weyl algebra $\W(U(\gg))$ generated by the enveloping algebra of a
given Lie algebra and the corresponding partial derivatives and give some example of these Weyl algebras.

Besides, we generalize  all considered objects (partial derivatives, Weyl algebra, de Rham complex) to the  Reflection Equation  algebra
and its modified version. This algebra and all related objects are called {\em braided} since they arise from braidings (see section 4). For
a more precise meaning of this term the reader is referred to \cite{GS3}. The explicit construction of the objects mentioned above is one
of the  purposes of the present paper.

From the other side, certain deformations of the usual derivative are known for a long time, for instance, the q-derivative (also called the
Jackson derivative) defined by
\be
\paq (f(t))=\frac{f(qt)-f(t)}{t(q-1)}
\label{qder}
\ee
or the difference operator
\be
\pah (f(t))=\frac{f(t+\h)-f(t)}{\h}
\label{hder}
\ee
(called below the $\h$-derivative) or their slight modifications. It is tempting to use them in order to introduce analogs of algebras whose
construction is based on the usual derivative. The most known examples are the q-Witt and q-Virasoro algebras.

The other purpose of the paper is to study deformation property of the enveloping algebra of the q-Witt algebra. We show that contrary to
the claim  of \cite{H}, the PBW property fails in this enveloping algebra. We demonstrate it in the section 5, as well as a similar claim for the
enveloping algebra of the $\h$-Witt algebra, constructed with the use of the $\h$-derivative instead of the usual one. Our reasoning is based
on the paper \cite{PP} where a version of the Jacobi condition useful for dealing with quadratic algebras and their quadratic-linear deformations
is exhibited\footnote{A version of this construction covering quadratic-linear-constant deformations of quadratic algebras was considered in
\cite{BG}.}. This condition is necessary for the PBW property and since it is not satisfied for the aforementioned  algebras,
we arrive to our conclusion.

In this connection we discuss other forms of the Jacobi condition which are useful for generalizing some other objects and operators associated with Lie
algebras, namely, the Chevalley-Eilenberg complex and the adjoint representation. It is worth noticing that in general these forms of the Jacobi
condition are not equivalent and each of them plays its own role in the theory of quadratic(-linear) algebras.

The paper is organized as follows. In the next section we compare different ways of associating a differential algebra with the enveloping algebra
of $gl(m)$.  In section 3 we discuss a generalization of this construction onto the enveloping algebras of  some other Lie algebras. In section 4 we
extend our construction to the Reflection Equation algebra. Here, the central problem consists in a convenient definition of the algebra generated by the
differentials of the generators of the initial algebra such that  the corresponding de Rham operator $d$ meets the usual property $d^2=0$. In section
5 we consider the aforementioned versions of the Witt algebra and show that the PBW property fails in their enveloping algebras. We complete the
paper (section 6) with a discussion on different forms of the Jacobi condition related to different generalizations of the Lie algebra notion.

{\bf Acknowledgement.}
The work of P.S. was supported by The National Research University Higher School of Economics' Academic Fund Program in 2014--2015, research
grant No 14-01-0173.

\section{Partial derivatives on $U(gl(m))$: different approaches}
\label{sec:two}

In what follows we deal with different deformations of the symmetric algebra $\Sym(\gg)$, where $\gg$ is a Lie algebra.
 Our main example is $\gg=gl(m)_\h$, where the subscribe $\h$ means that the parameter $\h$ is introduced as a multiplier in the $gl(m)$ Lie bracket. As usual, we fix a basis
$\{n_i^j\}$, $1\leq i,j\leq m$, in $gl(m)$ and the Lie brackets of the Lie algebra  $gl(m)_\h$   read
$$
[n_i^j,n_k^l]=\h (n_i^l \de_k^j-n_k^j \de_i^j),\qquad 1\leq i,j,k,l, \leq m.
$$

In each homogenous component of the algebra  $\Sym(gl(m))$ we fix a basis consisting of  symmetric elements, i.e. those invariant with respect to the
action of the symmetric group. Denote $\{e_\beta\}$ the corresponding basis of the whole algebra $\Sym(gl(m))$. Any element $e_\beta$ is a polynomial
in the generators of the algebra $\Sym(gl(m))$. A similar basis in the filtered quadratic-linear algebra $U(gl(m)_\h)$ will be denoted $\{\hat{e}_\beta\}$.
The element $\hat{e}_\beta$ can be obtained by replacing
the generators of $\Sym(gl(m))$ in the polynomial $e_\beta$
by  the corresponding generators of the algebra $U(gl(m)_\h)$.

Now, consider a linear map
$$
\al:\, \Sym(gl(m))\to U(gl(m)_\h)
$$
defined on the above bases as follows
$$
\al(e_\beta)= \hat{e}_\beta.
$$
This map is the central ingredient of the Weyl quantization method.

Using this map we can pull forward any operator ${\cal Q}:\, \Sym(gl(m))\to \Sym(gl(m))$ to that ${\cal Q}_\al:\,U(gl(m)_\h)\to U(gl(m)_\h)$ as follows
$$
{\cal Q}_\al=\al\circ {\cal Q}\circ\al^{-1}.
$$
In particular, we can pull forward partial derivatives from the algebra $\Sym(gl(m))$ to that $U(gl(m)_\h)$ and consider them as an appropriate noncommutative
analog of the usual partial derivatives. And visa versa, any operator defined in the algebra $U(gl(m)_\h)$ can be pulled back to $\Sym(gl(m))$. For instance,
the product in the latter algebra  being pulled back to the former algebra is called  $\star$-product (induced from
$U(gl(m)_\h)$. This product is often used in a quantization of dynamical
models (see \cite{K}). In such models (for example, the Schrodinger one) the kinetic part composed of momenta is classical but the usual product of coordinate
functions is replaced by the $\star$-product.

Equivalently, these models can be treated in terms of the algebra  $U(gl(m)_\h)$ but then the partial derivatives (momenta) of the kinetic part should be replaced
by their images with respect to the map $\al$. In fact, this method of defining the partial derivatives on the algebra   $U(gl(m)_\h)$ can be presented as follows.
One realizes a given element of the algebra $U(gl(m)_\h)$ in a symmetric form and employs the usual Leibniz rule to this element.

Note that other methods of defining the map $\al$ (for instance, the Wick one) also can be used. However, the way exhibited above (in fact, the Weyl quantization
method) gives rise to a $GL(m)$-covariant map.

 The second method of defining the partial derivatives on the algebra $U(gl(m)_\h)$ consists in modifying the Leibnitz rule. The modified Leibniz rule can be realized
via the coproduct defined  on the partial derivatives as follows\footnote{This form of the Leibniz rule was found by S.Meljanac and Z.{\v{S}}koda.}
\be
\De(\pa_i^j)=\pa_i^j\ot 1+1\ot \pa_i^j+\h\sum_k\pa^j_k\ot \pa_i^k.
\label{copr}
\ee
Hereafter, we use the notation $\pa_i^j=\pa_{n_j^i}$ for the partial derivative in the element $n^i_j$. Thus, we set by definition that
\be
\pa_i^j(n_k^l)=\de_i^l\de_k^j,
\label{def:partial}
\ee
i.e. this action is nothing but the pairing of the dual bases $\{n_i^j\}$ and $\{\pa_i^j\}$. Besides, we naturally assume the derivatives to be linear operators killing
elements of the ground field $\K$ which is assumed to be $\C$ or $\R$ depending on the context. Then, using the coproduct (\ref{copr}), one can extend the action
of the derivatives on polynomials in the generators.

The second form  of the Leibniz rule suggested in \cite{GS3} consists in the following. Consider an associative product $n_i^j\circ n_k^l=\de_k^j n_i^l$  in the
Lie algebra $gl(m)$. Note that $[n_i^j, n_k^l]=n_i^j\circ n_k^l- n_k^l\circ n_i^j$. Then, in addition to (\ref{def:partial}), the action of a derivative on a quadratic monomial in
generators is defined as follows:
$$
\partial_i^j(n_a^b\,n_c^d) = \partial_i^j(n_a^b)\,n_c^d+n_a^b\,\partial_i^j(n_c^d) +\hbar\,\partial_i^j(n_a^b\circ n_c^d).
$$
In general, the action of a derivative $\pa_i^j$ on a $p$-th order monomial $n_{i_1}^{j_1}...n_{i_p}^{j_p}$ gives rise to a sum of monomials whose order varies from
zero (a constant term) to $p-1$. In this sum the $(p-k)$-th order component $(1\le k \leq p)$ is composed from all monomials,
each of them  being obtained by the pairing of $\pa_i^j$ and the $\circ$-product of a
subset of $k$ elements from the initial monomial. Besides, the sum of all such   $(p-k)$-th order monomials has a multiplier $\h^{k-1}$. We illustrate this rule by an example of a third
order monomial:
\begin{eqnarray*}
\pa_i^j(n_a^b\, n_c^d\, n_k^l) \!\!\!&=&\!\!\! \pa_i^j(n_a^b)\, n_c^d\, n_k^l+n_a^b\, \pa_i^j(n_c^d)\, n_k^l+n_a^b\, n_c^d\, \pa_i^j(n_k^l)\\
&+&\!\!\!\h\left(\pa_i^j(n_a^b\circ n_c^d)\, n_k^l+\pa_i^j(n_a^b\circ n_k^l)\, n_c^d+n_a^b\,\partial_i^j (n_c^d\circ n_k^l)\right)+\h^2\pa_i^j(n_a^b\circ n_c^d\circ n_k^l).
\end{eqnarray*}

Observe that the partial derivatives commute with each other. Denote $\D$ the unital algebra generated by the partial derivatives. It becomes a bi-algebra being
equipped with the coproduct defined on the generators by formula (\ref{copr}) and the counit $\varepsilon : \D\to \K$ defined in the usual way: it kills all generators
$\pa_i^j$ and maps   $1_\D$ (the unit of $\D$) into the unit of the field.

The above coproduct allows one to introduce the so-called {\em permutation relations} between the partial derivatives and elements of the algebra ${\cal U}=U(gl(m)_\h)$
by the following rule
$$
\pa_i^j\ot n_k^l= (\pa_i^j)_1 \triangleright n_k^l\ot (\pa_i^j)_2,\quad {\rm where}\quad \De(\pa_i^j)=(\pa_i^j)_1\ot (\pa_i^j)_2
$$
in Sweedler's notation. Also, the notation $\triangleright$ stands for the action of an operator on an element.
Explicitly these permutation relations read:
$$
\pa_i^j\ot n_k^l-n_k^l\ot \pa_i^j=\de_i^l \de_k^j \,1_{\cal U} \ot 1_\D+ \h \, 1_{\cal U}\ot
(\pa_i^l \de_k^j-\pa_k^j \de_i^l).
$$
These permutation relations  can be presented in a matrix form as follows
\be
D_1\,P\, N_1\, P-P\, N_1\, P\, D_1=P+\h(D_1\,P-P\,D_1).
\label{perm}
\ee
Here $D=\|\pa_i^j\|$ and $N=\|n_i^j\|$ are the matrices composed of the elements $\pa_i^j$ and $n_i^j$ respectively (the low index labels the lines) and $A_1=A\ot 1$
for any matrix $A$. Also, $P$ stands for the matrix of the usual flip. Besides,  we omit the factors $1_{\cal U}$,  $1_\D$ and the sign $\ot$.

The algebra generated by two subalgebras $U(gl(m)_\h)$ and $\D$, equipped with the permutation relations  (\ref{perm}), is called the {\em Weyl algebra} and is denoted
$\W(U(gl(m)_\h)$. Note that for $\h=0$ we get the usual Weyl algebra generated by $\Sym(gl(m))$ and  the usual partial derivatives in the generators\footnote{Physicists
prefer to call this algebra the Heisenberg one.}.

The above permutation relations have been obtained via a passage to a limit $q\to 1$ in the permutation relations for the modified Reflection Equation  algebra under
assumption that the Hecke symmetry is a deformation of the usual flip. In general, the permutation relations themselves can be used for introducing partial derivatives. In
order to define the action of a derivative $\pa_i^j$ on an element $a\in U(gl(m)_\h)$ one   proceeds as follows.  One permutes the factors in the product $\pa_i^j\ot a$ by
means of the permutation relations and applies the counit to the right factor of the final element belonging to the tensor product $U(gl(m)_\h)\ot \D$.

Concluding this section, we resume that there are three ways of defining the partial derivatives on the algebra $U(gl(m)_\h)$.
One of them is based on using the coproduct (\ref{copr}), another one uses the product $\circ$ in the algebra $U(gl(m)_\h)$.
The third way is based on the related permutation relations.
 Similar ways also exist on the enveloping
algebras of the Lie super-algebras $gl(m|n)_\h$ and their "braided" analogs related to involutive braidings (see \cite{GS3}).
In section 4 we consider similar algebras related to non-involutive (namely, Hecke type) braidings. For them the only  way
based on permutation relations is known.

\section{Partial derivatives on other enveloping algebras}

The methods of defining partial derivatives on the algebra $U(gl(m)_\h)$ should be modified for other enveloping algebras. Consider some examples.

Let $\gh$ be a Lie subalgebra of the Lie algebra $gl(m)_\h$. In general, the above method of defining the partial derivatives via permutation relations fails since in the
permutation relations some extra-terms appear which do not belong to the subalgebra $\gh$. The same is true for the coproduct (\ref{copr}). Nevertheless, the partial
derivatives in elements of  $\gh$ are well defined  as operators. To show this, we fix a complementary subspace $W$ to $\gh$ that is $gl(m)_\h=\gh\oplus W$ as vector
spaces. Then we chose a basis in the  subalgebra $\gh$ and extend it to the basis in $gl(m)_\h$ by fixing a basis  in $W$. This new basis of the space $gl(m)$ is
subordinate to the direct sum $\gh\oplus W$. In the space generated by the partial derivatives we pass to the dual basis.

Let $x_i$ be an element of the chosen basis of $\gh$. Then, on applying the  derivative $\pa^{x_i}$ to a monomial composed of elements from $\gh$  we get a polynomial
also possessing this property though in the coproduct (\ref{copr}) presented in the new basis of $gl(m)$ certain external derivatives (i.e. derivatives in elements from $W$)
enter.

Now, consider some subalgebras of the Lie algebra $gl(2)_\h$ or more precisely, of its compact form $u(2)_\h$. On fixing in the latter algebra the standard basis $\{t, x,y,z\}$
such that
$$
[x,y]=\h z,\quad [y,z]=\h x,\quad [z,x]=\h y,\quad [t,x]= [t,y]=[t,z]=0,
$$
and taking in the dual space the basis $\pa^t$, $\pa^x$, $\pa^y$, $\pa^z$, we get the following permutation relations
$$
\begin{array}{l@{\quad}l@{\quad}l@{\quad}l}
 \pa^t\,t - t\, \pa^t =1+ {\h\over 2}\,\pa^t & \pa^t\, x - x\, \pa^t
=-{\h\over 2}\,\pa^x &
 \pa^t\, y - y\,  \pa^t=-{\h\over 2}\,\pa^y & \pa^t\, z - z\, \pa^t=- {\h\over 2}\,\pa^z\\
\rule{0pt}{7mm}
\pa^x\, t - t\,\pa^x = {\h\over 2}\,\pa^x &\pa^x \,x -  x\,\pa^x =1+ {\h\over 2}\, \pa^t &
\pa^x \, y-  y\,\pa^x = {\h\over 2}\,\pa^z & \pa^x \,z - z\, \pa^x  = - {\h\over 2}\,\pa^y \\
\rule{0pt}{7mm}
\pa^y \,t - t \, \pa^y = {\h\over 2}\,\pa^y & \pa^y \,x -  x\,  \pa^y = -{\h\over 2}\,\pa^z &
\pa^y \,y - y \,  \pa^y =1+ {\h\over 2}\, \pa^t & \pa^y \,z - z \,  \pa^ y= {\h\over 2}\,\pa^x\\
\rule{0pt}{7mm}
\pa^z \,t - t \,\pa^z = {\h\over 2}\,\pa^z & \pa^z \,x - x \,\pa^z = {\h\over 2}\,\pa^y&
\pa^z \,y -  y\,\pa^z = -{\h\over 2}\,\pa^x & \pa^z \,z - z \,\pa^z =1+ {\h\over 2}\, \pa^t.
\end{array}
$$

First, consider the Lie subalgebra $su(2)_\h \subset u(2)_\h$. As follows from \cite{GPS2} for any polynomial  of the form $f(x, y, z)=f_1(x)\,f_2(y)\,f_3(z)$
the action of the partial derivative $\pa^x$ is defined by
$$
\pa^x(f)=2\h^{-1} \,(B(f_1)\,A(f_2)\,A(f_3)+A(f_1)\,B(f_2)\,B(f_3)),
$$
where
$$
A(f(v))=\frac{1}{2}\Big(f(v-i\,\h/2)+f(v+i\,\h/2)\Big),\qquad B(f(v))=\frac{i}{2}\Big(f(v-i\,\h/2) - f(v+i\,\h/2)\Big).
$$
 Note that though the quantity $i=\sqrt{-1}$ enters these formulae, the result is real provided $f$ has real coefficients and $\h\in \R$. Similar formulae
are valid  for the derivatives $\pa^y, \pa^z$.  Consequently, we have
$$
\pa^x(f),\, \pa^y(f),\, \pa^z(f) \in U(su(2)_\h).
$$

\begin{definition} \rm Let $\gg$ be a Lie algebra and $U(\gg)$ its enveloping algebra. Choose a basis $\{x_i\}$, $1\leq i \leq m$ in $\gg$, and denote by
$c_{i,j}^k$ the structure constants of $\gg$ in this basis. Introduce an algebra $\W(U(\gg))$ generated by $U(\gg)$ and a commutative algebra $\D$ with
generators $\pa^l$, $1\leq l \leq m$, subject to the permutation relations
\be
[\pa^i,x_j]= - [x_j,\partial^i] =  b^i_{j,k}\, \pa^k+\de_j^i.
\label{br}
\ee
We call the algebra $\W(U(\gg))$ the Weyl algebra if  the Jacobi identity is valid for the bracket
$$
[\,,\,]: \wedge^2 (W)\to W\oplus  \K,\quad W=\gg\oplus \span(\pa^i),
$$
where the above bracket is defined by the initial Lie bracket on $\gg$, by the trivial bracket on $\D$ and by the bracket (\ref{br}) on $\span(\pa^i)\ot \gg$
and on $\gg \ot \span(\pa^i)$.
\end{definition}

Note that the Jacobi identity must be adapted to the case when the image of the bracket belongs to $W\oplus  \K$. In fact, we have only to satisfy the
relations
$$
[\pa^p, [x_i, x_j]]=[[\pa^p,x_i],x_j]-[[\pa^p,x_j],x_i],
$$
or in terms of the structure constants
$$
c_{i,j}^k b_{k,l}^p=b_{i,k}^p b_{j,l}^k-b_{j,k}^p b_{i,l}^k,\quad c_{i,j}^p=b_{i,j}^p-b_{j,i}^p.
$$

Observe that in virtue of the  PBW theorem the graded algebra $Gr\, \W(U(\gg))$ is canonically isomorphic to the commutative algebra generated by the
elements $x_i$ and $\pa^j$.

It is straightforward checking that the algebra $\W(U(u(2)_\h))$  is a Weyl algebra in the sense of the above definition. However, an attempt to define a similar
Weyl algebra for the enveloping algebra $U(su(2)_\h)$  as a quotient of $\W(U(u(2)_\h))$ fails. Indeed, it suffices to check that the relation
$$
[\pa^x,[x,y]]=[[\pa^x, x], y]-[[\pa^x,y],y]
$$
fails, if we assume that $\pa^t=0$. This example shows that in order to define the partial derivatives on an enveloping algebra we have, in general, to consider
the Weyl algebra related to a larger Lie algebra.

It is not the case for the subalgebra  $\gg\subset u(2)_\h$ generated by the elements $t$ and $x$. This Lie algebra is commutative. So, its enveloping algebra
coincides with $\Sym(\gg)$ and the corresponding partial derivatives and the Weyl algebra  can be defined in  the classical way. However, considering the subalgebra
of $\W(U(u(2)_\h))$ generated  by the elements $t,x,\pa^t,\pa^x$ we get another Weyl algebra corresponding to the same Lie algebra $\gg$. Thus, we get  two
different Weyl algebras related to the algebra $\gg$.

In general, for the two-dimensional commutative Lie algebra the permutation relations must be of the form
$$
\pa^t t-t\pa^t=1+(a_1\pa^t+b_1\pa^x),\quad \pa^t x-x\pa^t= (a_2\pa^t+b_2\pa^x),
$$
$$
\pa^x t-t\pa^x=(c_1\pa^t+d_1\pa^x),\quad\pa^x x-x\pa^x=1+ (c_2\pa^t+d_2\pa^x).
$$

It would be interesting to classify all possible families of constants $a_1,\dots ,d_2$ giving rise to the Weyl algebras on the two-dimensional commutative Lie algebra
(as well as on two-dimensional noncommutative one). Two examples above correspond to the following families of the constants. In the classical case all constants
are trivial. In the other one nontrivial constants are: $a_1=d_1=-b_2=c_2=\frac{\hh}{2}$.

We get a little bit more general Weyl algebra by putting
$$
\qquad b_1=c_1=a_2=0,\quad a_1=d_1=c_2
$$
(here we do not impose an restriction on $b_2$).

In conclusion, we want to mention the following fact. Though the way of introducing the partial derivatives on an enveloping algebra via permutation relations is not universal,
it is more general than the way based on the coproduct defined on the algebra $\D$. This observation is also valid for braided algebras considered in the next section.
We have not succeeded in finding a coproduct corresponding to the  permutation relations (\ref{permm}).

\section{Braided Weyl algebras and related de Rham complex}

In this section we consider a {\em braided analog} of the enveloping algebra $U(gl(m)_\h)$ and differential calculus on it. By braided analog we mean the so-called
Reflection Equation (RE) algebra in its modified form. Let us recall the definition of this algebra.

Let $V$ be a vector space over the ground field $\K$ and $R:\vv\to\vv$ be a linear invertible operator satisfying the braid relation
$$
(R\ot I)(I\ot R)(R\ot I)=(I\ot R)(R\ot I)(I\ot R)
$$
(in an equivalent form it is also called the {\em quantum Yang-Baxter equation}). Such an operator $R$ is called a {\em braiding}. If a braiding $R$ is subject to an
additional condition
$$
(R-q\, I)(R+\qq\, I)=0,\quad q\in \K,
$$
it is called a {\em Hecke symmetry} provided $q\not=1$ and an {\em involutive  symmetry} provided  $q=1$.

By modified Reflection Equation  algebra we mean a unital algebra generated by elements $n_i^j$, $1\leq i, j \leq m$, subject to the system of relations
\be
R\, N_1 R\, N_1-N_1 R\,N_1R=\h\,(R\, N_1-N_1 R),\quad \h\in \K
\label{mRE}
\ee
where $N=\|n_i^j\|$ and $N_1=N\ot I$. We omit the term "modified" if $\h=0$. The algebra (\ref{mRE}) will be denoted $\N(q,\h)$ provided $\hh\not=0$ or $\N(q)$
provided $\h=0$.

Below we assume $R$ to be a {\em skew-invertible} Hecke symmetry. This means that there exists an operator $\Psi:\vv\to\vv$ such that
\be
{\rm Tr}_2 R_{12}\Psi_{23}= {\rm Tr}_2 \Psi_{12}R_{23}=P_{13},
\label{def:psi}
\ee
where $\rm Tr$ stands for the usual trace, and indices label the spaces where the operators act. In what follows we shall need the operators
\be
B={\rm Tr}_1\Psi_{12},\qquad C = {\rm Tr}_2\Psi_{12}.
\label{def:B}
\ee
As a direct consequence of the definition of the operator $\Psi$ we have
\be
{\rm Tr}_1B_1R_{12} = I,\qquad {\rm Tr}_2 C_2R_{12} = I.
\label{BCR}
\ee
The operator $B$ is supposed to be invertible. Then, it can be shown that
\be
B\cdot C =  q^{-2m}I, \quad {\rm Tr}B = {\rm Tr}C = q^{-m} m_q, \quad m_q=\frac{q^m-q^{-m}}{q-q^{-1}}.
\label{BC-prop}
\ee

Note, that the algebra (\ref{mRE}) is filtered. We call it the quadratic-linear one since it is defined by the quadratic-linear relations. We treat this algebra as a {\em braided
analog} of the enveloping algebra $U(gl(m)_\h)$. Indeed, it is possible to define a {\em braided Lie bracket} such that the modified RE algebra has the sense of the
enveloping algebra of the corresponding braided Lie algebra (see \cite{GS1}). Furthermore, the algebra  $\N(q,\h)$ has the following  properties (see \cite{GPS1}):

1. If $R$ comes from the quantum group (QG) $U_q(sl(m))$, the algebra  $\N(q,\h)$ has a good deformation property. This means that the homogeneous components
$\N^{k}(q)$, $k=0,1,2,\dots $, of the  algebra $\N(q)$ have the classical dimensions, i.e.
$$
\dim\, \N^{k}(q)=\dim U^{k}(gl(m)_\h)
$$
for any $k$ and a generic $q$. Also, for the algebra $\N(q,\h)$ there is a sort of the PBW theorem ensuring that the associated graded algebra $Gr (\N(q,\h))$ is isomorphic
to $\N(q)$ (for a discussion on the PBW property see the next section).

2. This algebra can be equipped with a braided bi-algebra structure (see \cite{GPS2} for a definition). This structure is determined by the usual counit and the coproduct
such that for $\h=1$ it has the form
$$
\De(n_i^j)= n_i^j\ot 1+ 1\ot n_i^j-(q-\qq)  \sum _k n_i^k\ot n_k^j.
$$
But similarly to super-algebras, the product of two such elements $\De(n_i^j)\De(n_k^l)$ contains an operator transposing two middle factors in the product. The transposing
operator depends on a concrete Hecke symmetry, defining the algebraic structure of  $\N(q,\h)$.

3. The representation of this algebra is similar to that of $U(gl(m))$ or  $U(gl(m|n))$ (depending on $R$), an analog of the adjoint representation included.

4. The structure of the center of $\N(q,\h)$ is similar to that of $U(gl(m))$ (or $U(gl(m|n))$).

5. If $R$ comes from $U_q(sl(m))$, the algebra $\N(q,\h)$ is covariant with respect to the action of this QG. In general, a similar property can be formulated via a coaction
of the RTT algebra.

Note that if $R$ is an involutive symmetry, the algebra $\N(q,\h)$ turns into the enveloping algebra of a generalized Lie algebra, introduced by one of the authors in the 80's.

The point is that a braided analog of the partial derivatives can be introduced on the algebra $\N(q,\h)$ (first it was done in \cite{GPS2}). This enable us to define a braided
analog of the Weyl algebra corresponding to the algebra  $\N(q,\h)$.

\begin{definition}
\label{def:weyl-alg}
\rm  {\em The braided Weyl algebra} $\W(\N(q,\h))$ is an associative unital algebra generated by two subalgebras $\N(q,\h)$ and $\D$ provided that the
following conditions are satisfied:
\begin{enumerate}
\item
As a vector space the algebra $\W(\N(q,\h))$ is isomorphic to $\N(q,\h)\ot \D$.
\item
The subalgebra $\D$ is generated by elements $\pa_i^j$, $1\le i,j\le m$, subject to the following relations
\be
R^{-1} D_1 R^{-1} D_1-D_1 R^{-1} D_1 R^{-1}=0,
\label{ddd}
\ee
where $D=\|\pa_i^j\|$ and $D_1=D\ot I$.
\item
The permutation relations between the generators $n_i^j$, $1\le i,j\le m,$ of the subalgebra $\N(q,\h)$ and the generators $\pa_i^j$ of the subalgebra $\D$ are as follows
\be
D_1 R\,N_1 R-R \, N_1 R^{-1} D_1=R+\h D_1 R.
\label{permm}
\ee
\end{enumerate}
\end{definition}

In the limit $q=1$ (provided that $R$ is a deformation of the usual flip $P$) the relations (\ref{ddd}) turn into the commutativity conditions for the generators $\partial_i^j$
and the equalities (\ref{permm}) turn into the permutation relations (\ref{perm}). Note that the relations (\ref{ddd}) and (\ref{permm}) have been introduced  in \cite{GPS2}.

Now, we have to give an operator meaning to the elements of the subalgebra ${\cal D}$ since we intend to interpret the generators $\partial_i^j$ as analogs of partial
derivatives. The permutation relations (\ref{permm}) in the above Definition \ref{def:weyl-alg} allows one to define an action of the subalgebra ${\cal D}$ on the subalgebra
${\cal N}(q,\h)$. On the level of generators this action is as follows:
\be
\partial_i^j(n_k^p) = \delta_i^p B_k^j,
\label{action-der}
\ee
where $\|B_i^j\|$ is the matrix of the operator $B$ introduced in (\ref{def:B}). On an arbitrary monomial in generators $n_i^j$ the action is extended with the help of the
permutation relations (\ref{permm}) or by means of the same scheme as in the section \ref{sec:two} with the use of the counit defined on the algebra $\D$.  Let us point out
that this action gives a representation of the subalgebra ${\cal D}$ in the {\it algebra} $\N(q,\h)$, that is  the action (\ref{action-der}), extended on the whole algebra ${\cal N}(q,\h)$ via (\ref{permm}),    respects the algebraic structure of
$\N(q,\h)$. In the classical limit $q=1$ the operator $\pa_i^j$ turns into the usual partial derivative in  $n_j^i$: $\pa_i^j=\partial/\partial n_j^i$.

Our next aim is to define the space of differential forms on the algebra   $\N(q,\h)$ and to introduce an analog of the de Rham operator on it. First of all, we need a
braided analog of the external algebra generated by the differentials $dn_i^j$. In the classical case it is identified with the skew-symmetric algebra $\bigwedge(gl(m))$.
In our current setting we define the corresponding analog $\bigwedge_q$ as a quotient of the free tensor algebra generated by the linear space $\span(d n_i^j)$ over
the ideal generated by the matrix elements of the left hand side of the equality
\be
R_{12} \Om_1 \hat{\Psi}_{12} \Om_1+\Om_1 \hat{\Psi}_{12} \Om_1 R_{12}^{-1} = 0.
\label{skew}
\ee
Here as usual $\Om_1=\Om\ot I$ and the matrix elements $\Omega_i^j$ of the $m\times m$ matrix $\Omega  = \|\Omega_i^j\|$ are the linear combinations of the
differentials:
\be
\Omega_i^j = (B^{-1})_i^k\,dn_k^j,
\label{omega-matr}
\ee
where the summation over the repeated index is understood. The symbol $\hat{\Psi}$ stands for the following operator
$$
\hat\Psi_{12} = \Psi_{21} +(q-q^{-1})q^{2m}B_1C_2.
$$
Using the definitions and properties (\ref{def:psi})--(\ref{BC-prop}) one can easily verify that
\be
{\rm Tr}_1 \hat{\Psi}_{12} R^{-1}_{13}=P_{23} = {\rm Tr}_1\hat \Psi_{21}R^{-1}_{31}.
\label{psi-hat}
\ee
Below, we give a motivation for the definition (\ref{skew}) of the algebra $\bigwedge_q$.

Consider the product ${\frak D}_q=\bigwedge_q\ot \N(q,\h)$ which is a right $\N(q,\h)$-module. In order to convert this module into an associative algebra,
we have to introduce some permutation relations between the algebras $\bigwedge_q$ and $\N(q,\h)$. However, we do not use this structure and shall consider the
space ${\frak D}_q$ as a one-sided $\N(q,\h)$-module only. Elements of this module are called  {\em braided differentials}. Elements of the right $\N(q,\h)$-module
${\frak D}_q^{k}=\bigwedge_q^k\ot \N(q,\h)$ are called  {\em braided $k$-differentials}. Here, as usual, $\bigwedge_q^k$ stands for the $k$-th degree homogenous
component of the quadratic algebra $\bigwedge_q$.

Now, define the braided analog of the de Rham operator $d: {\frak D}_q^k\rightarrow {\frak D}_q^{k+1}$. Let
$$
\om= \om_0 \ot f,\quad f\in \N(q,\h),\,\, \, \om_0\in {\bigwedge}_q^k
$$
be a $k$-differential: $\omega\in {\frak D}_q^k$. Then we set by definition:
\be
d\, \om= \om_0 \ot
 \sum_{i,j} \Om_i^j \ot \pa^i_j (f)\in {\frak D}_q^{k+1}.  
\label{def:d}
\ee

In fact, the map $d$ consists in inserting the element $\sum_{i,j}\Om_i^j \ot \pa^i_j$ inside of the $k$-form $\om$ with subsequent application of the partial derivatives to the
element $f\in \N(q,\h)$. The following claim is the main motivation of our definition (\ref{skew}) of the algebra $\bigwedge_q$.

\begin{proposition}
The usual property  $d^2=0$ holds.
\end{proposition}

\noindent
{\bf Proof.} First, let us recall some facts from the theory of monoidal categories. Let ${\frak A}$ be a monoidal rigid category of finite dimensional vector spaces and $U$ be
its object. Let $U^*$ be its right dual (see  \cite{CP} for detail). This means that there exists an {\it evaluation  map} $U\ot U^*\to \K$ and a {\it coevaluation map} $\K\to U^*\ot U$
which are in a sense coordinated. Let $\{u_i\}$ be a basis of $U$ and $\{u^j\}$ be its right dual, i.e. $\langle u_i, u^j \rangle= \de_i^j$. Then as follows from the definition, the
coevaluation map is generated by $1\mapsto \sum_k  u^k\ot u_k$.

Now, consider a subspace $I\subset U^{\ot 2}$ and the quadratic algebra $\Sym(U)=T(U)/\langle I \rangle$ playing the role of the symmetric algebra of the space $U$. Also,
consider the subspace $I^\bot \subset (U^*)^{\ot 2}$ orthogonal to $I$ with respect to the pairing
\be
\langle x\ot y, z\ot v \rangle=\langle x, v\rangle\,\langle y, z \rangle,\quad x, y\in U,\,\, z, v \in U^*.
\label{ppa}
\ee
The algebra $\bigwedge(U^*)=T(U^*)/\langle I^\bot \rangle$ plays the role of the skew-symmetric algebra of the space $U^*$.

Let us form a complex
$$
\delta:{\bigwedge}^k(U^*)\ot \Sym^l(U) \to {\bigwedge}^{k+1}(U^*)\ot \Sym^{l+1}(U),
$$
where the map $\delta$ is defined by
$$
u^{j_1}\ot\dots\ot  u^{j_k}\ot u_{i_1}\ot \dots\ot u_{i_l}
\stackrel{\delta}{\mapsto} u^{j_1}\ot\dots\ot  u^{j_k}\ot  \sum_m (u^m \ot u_m)\ot u_{i_1}\ot \dots\ot u_{i_l}.
$$

Let us emphasize that the map $\delta$ consists in introducing the unit 1 inside of the product $\bigwedge^k(U^*)\ot\Sym^l(U)$ with subsequent applying  the coevaluation
operator to the unit.

\begin{lemma}
The following property holds: $\delta^2=0$.
\end{lemma}

\noindent
{\bf Proof.} In order to prove the lemma, we have  to show that the element
\be
\sum_{m,n} u^m\ot u^n \ot u_n\ot u_m,
\label{elem}
\ee
which corresponds to the operator $\delta$ applied twice,
vanishes in the product $\bigwedge^2(U^*)\ot \Sym^2(U)$.

Let $\{z_i\}$ be a basis of the subspace $I\subset U^{\ot 2}$ and $\{\overline{z^j}\}$ be a basis of the subspace $I^\bot\subset (U^*)^{\ot 2}$. Let us complete the former basis
up to a basis $\{z_i, \overline{z_j} \}$ of the whole space $U^{\ot 2}$ and the latter one up to a basis $\{{z^i}, \overline{z^j}\}$ of the space $(U^*)^{\ot 2}$ so that the basis
$\{z^i, \overline{z^j}\}$ be the right dual of that $\{z_i, \overline{z_j} \}$. Then the element (\ref{elem}) can be presented as follows
$$
\sum_i z^i\ot  {z_i} +\sum_j \overline{z^j}\ot \overline{z_j}.
$$
It is clear that this element vanishes in the product $\bigwedge^2(U^*)\ot \Sym^2(U)$.

\medskip

Now, go back to the proposition. We treat the space $\span(\pa_i^j)$ as an object $U$ from the above lemma and the space $\span(dn_i^j)$ as its right dual $U^*$. Besides,
the basis $\{\Omega_i^j\}$ (\ref{omega-matr}) of the space $U^*$ is the right dual to $\{\pa_k^l\}$ with respect to the pairing
\be
\langle D_1,\Omega_2\rangle = P_{12}\quad {\rm or}\quad \langle \pa_i^j, \Omega_k^p \rangle=\de_i^p\, \de_k^j.
\label{D-O}
\ee
In fact, if we identify $\span(n_k^l)$ and $\span(d n_k^l)$ as linear spaces, this pairing is nothing but the action of the partial derivatives on the generators of the RE algebra
in the spirit of the classical differential calculus.

The role of the subspace $I\subset U^{\ot 2}$ is played by the left hand side of (\ref{ddd}), giving rise to the RE algebra but with the braiding $R^{-1}$ instead of $R$. The only
claim has to be shown is that the left hand side of  (\ref{skew}) is just $I^\bot$. In order to prove this, we fix the basis $X_{12} = R^{-1} D_1 R^{-1} D_1$ in the space $U$ and that
$X^*_{12} = \Om_1 \hat{\Psi}  \Om_1 R$ in the space $U^*$.

\begin{lemma}
\label{lem:5}
The basis $X^*_{12}$ is right dual to $X_{12}$, that is $\langle X_{12}, X^*_{34}\rangle = P_{13}P_{24}$.
\end{lemma}

\noindent
{\bf Proof.} The claim of the lemma is verified by a direct calculation on the base of (\ref{psi-hat}), (\ref{ppa}) and (\ref{D-O}):
\begin{eqnarray*}
\langle X_{12}, X^*_{34}\rangle& = &\langle R_{12}^{-1}D_1R_{12}^{-1}\langle D_1\,,\,\Omega_3\rangle\hat\Psi_{34}\Omega_3R_{34}\rangle =
\langle R_{12}^{-1}D_1R_{12}^{-1}\,,\,\hat\Psi_{14}\Omega_1R_{14}\rangle P_{13}\\
& = & R_{12}^{-1}{\rm Tr}_0\Big( P_{01}R_{12}^{-1}\hat\Psi_{14}\langle D_0\,,\,\Omega_1\rangle\Big) R_{14} P_{13} =
R_{12}^{-1}{\rm Tr}_0\Big( P_{01}R_{12}^{-1}\hat\Psi_{14}P_{01}\Big) R_{14} P_{13} \\
&=& R_{12}^{-1}{\rm Tr}_0\Big( R_{02}^{-1}\hat\Psi_{04}\Big) R_{14} P_{13} = R_{12}^{-1}P_{24} R_{14} P_{13} = P_{24}P_{13}.
\end{eqnarray*}

Introduce now two operators $Q:U^{\ot 2}\to U^{\ot 2}$ and ${Q'}:U^{\ot 2}\to U^{\ot 2}$ defined as follows
$$
Q(R^{-1} D_1 R^{-1} D_1)=D_1 R^{-1} D_1R^{-1},\quad
{Q'}(R^{-1} D_1 R^{-1} D_1)=D_1 R^{-1} D_1 R.
$$
Similar operators were considered in \cite{GS2}. They were also used in the construction of a differential calculus on a q-Minkowski space algebra in \cite{M1, M2}.
It is easy to see that $Q$ and $Q'$ commute with each other and satisfy the following relation
$$
(I-Q)(I+Q')=0.
$$
Also, the subspace $I\subset U^{\ot 2}$ defined by (\ref{ddd}) can be written as follows
$$
(I-Q)(R^{-1} D_1 R^{-1} D_1)=0.
$$
Whereas, the equation
$$
(I+Q')(R^{-1} D_1 R^{-1} D_1)=0
$$
defines the subspace of $U^{\ot 2}$ which is complementary to $I$. In fact, the subspace $I^{\bot}$ is just $(I+{Q'}^*)(U^*)^{\ot 2}$ where ${Q'}^*$ stands for the
conjugate operator to $Q'$.

This shows that the spaces defined respectively by the left hand side of (\ref{ddd}) and (\ref{skew}) are orthogonal to each other\footnote{This can be also verified by a
direct calculation, similar to that in the proof of Lemma \ref{lem:5}.}. The fact that they are maximal (i.e. the latter space include all elements orthogonal to the former
space) can be shown from considering the dimensions of these spaces (first, for a generic $q$ with subsequent passage to all $q$).

This completes constructing an analog of the de Rham complex corresponding to the modified RE algebra.\hfill\rule{6.5pt}{6.5pt}

\begin{remark}\rm  Note that we precise no way of completing either the basis of the subspace $I$ or that of the space $I^\bot$. Nevertheless, there exists the "most
natural" choice to do so or, equivalently, to fix a complimentary subspace to $I$ in the space $U^{\ot 2}$. If the space $I\subset U^{\ot 2}$ is defined by the left hand side
of (\ref{mRE}) then we define its complementary subspace as
\be
R\, N_1 R\, N_1 + N_1 R\, N_1 R^{-1}.
\label{syme}
\ee
The quotient of the tensor algebra of $U$ over the ideal generated by this subspace is often treated to be a braided analog of the skew-symmetric algebra of $U$ (see
\cite{GS2}).
\end{remark}

\section{q-Witt algebra:  deformation property}

In this section we deal with the so-called q-Witt algebra ($q$ is assumed to be generic). This algebra is usually defined in the same way as the classical Witt algebra is but
with the q-derivative (\ref{qder}) instead of the usual one. Let us precise that $\pa_q(x^k)=k_q x^{k-1}$, $k\in \Z$. Hereafter, we use the notation $m_q=\frac{q^m-1}{q-1}$.

Note that the Leibniz rule for the q-derivative reads 
\be \pa_q(f(x) g(x))=(\pa_qf(x))g(x)+f(qx)\pa_qg(x), \label{modifL} \ee
whereas its permutation relation with the
generator $x$ is:
\be
\pa_q\, x -q\, x\, \pa_q=1.
\label{perr}
\ee
Below, we do not use the  Leibniz rule (\ref{modifL}) (see remark in the end of the section). 

Now, similarly to the usual Witt algebra, consider the  operators
$$
e_k=x^{k+1}\, \pa_q, \quad k\in {\Bbb Z}
$$
acting on the algebra $\K[x,x^{-1}]$. These operators act on the elements $x^l$ as follows
$$
e_k(x^l)= l_q\, x^{k+l},\quad  l\in {\Bbb Z}
$$
and are subject to the  relations
\be
q^{m+1} e_m e_n-q^{n+1} e_n e_m-((n+1)_q-(m+1)_q) e_{m+n}=0.
\label{rela}
\ee

These relations are usually considered (see \cite{H} and the references therein) as a motivation for introducing the following  "q-Lie bracket"
\be
U\ot U\to U:\,\, e_m\ot  e_n \mapsto [e_m, e_n]=((n+1)_q-(m+1)_q) e_{m+n},
\label{qLie}
\ee
where $U=\span(e_k)$ is the space of all finite linear combinations of the elements $e_k$. Then by q-Witt algebra one means the space $U$
endowed with the q-Lie bracket  (\ref{qLie}), which is assumed, of course to be a bilinear operator.  We denote this q-Witt algebra $\Wq$. Its
enveloping algebra $U(\Wq)$ is defined to be  the quotient of the free tensor algebra of the space $U$ over the ideal generated by the left hand
side of (\ref{rela}).

Emphasize that the bracket (\ref{qLie}) is well-defined on the whole space $U^{\ot 2}$. This bracket has the following properties:
\begin{enumerate}
\item The "q-skew-symmetry":
$$
[e_m, e_n]=-[e_n, e_m];
$$
\item The "q-Jacobi relation":
$$
(1+q^k) [e_k,[e_l, e_m]]+(1+q^l) [e_l,[e_m, e_k]]+(1+q^m) [e_m,[e_k, e_l]]=0.
$$
\end{enumerate}
The first relation entails that the element $e_m\ot  e_n + e_n\ot  e_m$ is killed by the bracket. Consequently, we have two subspaces in the space
$U^{\ot 2}$
\be
I_+=\span(e_m\ot  e_n + e_n\ot  e_m),\quad I=I_-=\span(q^{m+1} e_m e_n-q^{n+1}  e_n e_m),
\label{II}
\ee
which are analogs of symmetric and skew-symmetric subspaces (in fact, the symmetric one is classical).

Below, we deal with the PBW theorem in the form suggested in \cite{PP}. Namely, let  $U$ be a finite dimensional vector space over the field $\K$ and
$I\subset U^{\ot 2}$ be a subspace. Consider an operator $[\,\,,\,]: I\to U$ satisfying two conditions
\begin{enumerate}
\item  $[\,\,,\,]_{12}\otimes {\rm id}_3-{\rm id}_1\otimes [\,\,,\,]_{23}:\quad I\ot U \bigcap U\ot I\to I$;

\item $[\,\,,\,]\circ ([\,\,,\,]_{12}\otimes {\rm id}_3-{\rm id}_1\otimes [\,\,,\,]_{23}):\quad I\ot U \bigcap U\ot I\to 0,$
\end{enumerate}
where in the second line the symbol $\circ$ means the composition of the maps. (Below, we omit this symbol as well as
the identical operators.) 

If, in addition, the quadratic algebra $\A=T(U)/\langle I \rangle$
is Koszul then the associated graded algebra $Gr \A_{[\,,\,]}$ where $\A_{[\,,\,]}=T(U)/\langle I-[\,,\,] I\rangle$ is canonically isomorphic to $\A$. Here
$\langle I \rangle$ stands for the ideal generated by a set $I$ and by $I-[\,,\,] I$ we mean the family of elements $u-[\,,\,]u,\,\, u\in I$.

This is just the PBW theorem under the form of \cite{PP}. Below, the call the both conditions listed above the {\em Jacobi-PP condition}.

Emphasize that the subspace $I\ot U \bigcap U\ot I\subset U^{\ot 3}$ is an analog of the space of third degree skew-symmetric elements. Note, that the
bracket is defined only on the subspace $I$. Thus, the first of the above conditions (which means that the bracket maps $I\ot U \bigcap U\ot I$  into $I$),
ensures a possibility to apply the bracket once more.

Let us also show that the first condition above (without assuming the algebra  $T(U)/\langle I \rangle$ to be Koszul) is necessary for the canonical isomorphism.
Consider the  element
\be
([\,,\,]_{12} - [\,,\,]_{23}) Z,
\label{term}
\ee
where $Z$ is an arbitrary element belonging to $ I\ot U \bigcap U\ot I$. Since the element $Z-Z$ equals to 0 in the algebra  $\A_{[\,\,,\,]}$, its image under replacing
factors from $I\ot U$ (resp., $U\ot I$) by the terms  $[\,,\,]_{12} Z$ (resp., $[\,,\,]_{23} Z$) is also trivial in the algebra  $\A_{[\,,\,]}$. If nevertheless, the term
(\ref{term}) does not belong to $I$, we have that there is an element which is trivial in $Gr \A_{[\,,\,]}$ and is nontrivial in $\A$. Consequently, the canonical isomorphism
of the algebras $Gr \A_{[\,,\,]}$ and $\A$ does not exists.

\begin{remark} \label{Rem} \rm Note that to describe the space $I^{(3)}= I\ot U \bigcap U\ot I\subset U^{\ot 3}$ explicitly is not an easy deal in general. However,
if  the subspace $I\subset U^{\ot 2}$ is generated by elements of the form
\be
e_i e_j-c(i, j) e_j e_i,\quad c(i, j)\not=0\,\,\, \forall \, i, j
\label{easy}
\ee
the space $I^{(3)}$ is easy to describe. First, consider the case $\dim\, U=3$. Let $\{x,y,z\}$ be a basis of the space $U$. We set
$$
I=\span(xy-a yx,\, yz-b zy,\, zx-c xz),\quad a\,b\, c\not=0.
$$
Then the space $I^{(3)}$ is one-dimensional and is generated by the following element
\begin{eqnarray}
Z(x,y,z) &=&c(xy-a yx)z+a(yz-b zy)x+b(zx-c xz)y\nonumber \\
&=&b z(xy-a yx)+c x (yz-b zy)+a y (zx-c xz).
\label{Zet}
\end{eqnarray}
If $\dim\, U >3$, the space $I^{(3)}$ is generated by all elements  $Z(e_k, e_l, e_m)$, each of them being associated with a triple $e_k, e_l, e_m$.
\end{remark}

Now, go back to the q-Witt algebra. This algebra is infinite dimensional. However, if by $U$ we mean all finite linear combinations of the generators  $\{e_i\}$, and
by $U^{\ot k}$ we also mean the finite linear combinations of $e_{i_1}\ot e_{i_2}\ot\dots \ot e_{i_k}$,  then we can extend our reasoning to this case.

Namely, denote the vector space of finite linear combinations of elements $q^{k+1} e_k\, e_l- q^{l+1} e_l\, e_k$ by $I$ and consider an element belonging to the
space $I\ot U\bigcap U\ot I$
$$
Z=q^{l+1} q^{m+1} (q^{l+1} e_l e_m-q^{m+1} e_m e_l)e_k + q^{m+1} q^{k+1} (q_{m+1} e_m e_k-q^{k+1} e_k e_m)e_l+
$$
$$
q^{k+1} q^{l+1} (q_{k+1} e_k e_l-q^{l+1} e_l e_k)e_m = q^{2(m+1)} e_m(q^{k+1} e_k e_l-q^{l+1} e_l e_k) +
$$
$$
q^{2(k+1)} e_k (q^{l+1} e_l e_m-q^{m+1} e_m e_l) +q^{2(l+1)} e_l (q^{m+1} e_m e_k-q^{k+1} e_k e_m).
$$

Compute the images of this element under the maps $[\,,\,]_{12}$ and $[\,,\,]_{23}$ correspondingly. We have
$$
[\,,\,]_{12}Z=q^{l+1} q^{m+1} ((m+1)_q-(l+1)_q) e_{l+m}e_k+q^{m+1} q^{k+1}((k+1)_q-(m+1)_q) e_{m+k}e_l+
$$
$$
q^{k+1} q^{l+1}((l+1)_q-(k+1)_q) e_{k+l}e_m,
$$
$$
[\,,\,]_{23}Z=q^{2(m+1)}((l+1)_q-(k+1)_q) e_m  e_{k+l}+q^{2(k+1)} ((m+1)_q-(l+1)_q) e_k e_{l+m}+
$$
$$
q^{2(l+1)}(((k+1)_q-(m+1)_q) e_l  e_{m+k}.
$$

Let us assume that the numbers $k$, $l$, $m$, $k+l$, $k+m$ and $l+m$ are pairwise distinct. Then the difference $[\,,\,]_{12}Z-[\,,\,]_{23}Z$ belongs to $I$ iff it is so for
the element
\be
q^{l+1} q^{m+1} ((m+1)_q-(l+1)_q) e_{l+m}e_k-q^{2(k+1)} ((m+1)_q-(l+1)_q) e_k e_{l+m}
\label{coeff}
\ee
and for two similar elements obtained by cyclic permutations $k\to l \to m$. However, it is evident that for a generic $q$ the element (\ref{coeff}) does not belong to $I$ since
the vector 
$$(q^{l+1} q^{m+1} ((m+1)_q-(l+1)_q),\, -q^{2(k+1)} ((m+1)_q-(l+1)_q))$$
 composed of the coefficients of this element is not collinear to $(q^{l+m+1},\, -q^{k+1})$.

Thus, the first of the above conditions is not satisfied the algebra ${\rm Gr}(U(\Wq))$ is not isomorphic to $T(U)/\langle q^{k+1} e_k e_l-q^{l+1} e_l e_k \rangle$.
\begin{remark} \rm
If we introduce a parameter $\h$ as a multiplier in the right hand side of the bracket (\ref{qLie}) of the q-Witt algebra we get a two parametric analog of the usual Witt algebra.
On putting $\h=0$ we get a quadratic algebra $T(U)/\langle q^{k+1} e_k e_l-q^{l+1} e_l e_k \rangle$ which possesses a good deformation property. Since this quadratic algebra
is infinite dimensional, we should precise the meaning of this property. The ordered monomials
$$
e_1^{k_1}e_2^{k_2}\dots e_l^{k_l},\qquad k_1+k_2+\dots +k_l=k
$$
form a basis of its $k$-th degree homogeneous component. It can be considered as a quantization of the corresponding Poisson structure (the proof is left for the reader).
Nevertheless, the q-Witt algebra is not a two-parameter quantization of a Poisson pencil. This is due to the fact that the passage from the mentioned quadratic algebra to its
filtered (quadratic-linear) analog is not a deformation.
\end{remark}

In a similar way we can introduce another analog of the Witt algebra, called below $\h$-Witt one. In its construction the usual derivative is replaced by its difference analog
(\ref{hder}). The permutation relation with $x$ reads
$$
\pa_\h x-x \pa_\h=1+\h \pa_\h.
$$
Note that the algebra generated by $x$ and $\pa_\h$ is a Weyl algebra but it is not so for the algebra generated by $x$ and $\pa_q$. (However, it is a Weyl algebra in the sense of a more general definition exhibited in the next section.)

Now, consider operators $e_k= \exp(ikx)\pah$, $k\in {\Bbb Z}$ acting onto the space of real continuous functions. By using the permutation relation
$$
\pah\, \exp(ikx)-\exp(ik\h)\,\exp(ikx)\,\pah=\frac{\exp(ik\h)-1}{\h}\, \exp(ikx),
$$
we get the following relations between these operators
$$
\exp(ik\h)\,e_k \,e_l-\exp(il\h)\,e_l \,e_k=\frac{\exp(il\h)-\exp(ik\h)}{\h}\, e_{k+l}.
$$
So, denoting $q=\exp(i\h)$, we can see that the difference of this structure from that described above is unessential.

Also, it is tempting to introduce an $\h$-analog of the Lie bracket  by putting
$$
[e_k, e_l]=\frac{\exp(il\h)-\exp(ik\h)}{\h} e_{k+l}.
$$
Finally, for the reason presented above,  the PBW theorem in the enveloping algebra of this "generalized Lie algebra" fails. The detail is left to the reader.

Two above analogs of the usual derivatives being put together give rise to a $(q,\h)$-counterpart of the derivatives. This $(q,\h)$-derivative have the
following permutation relation with $x$:
$$
\pa_{q,\h} \,x -qx\, \pa_{q,\h}=1+\h \pa_{q,\h}.
$$
This permutation relation can be deduced from that for the q-derivative via the change of the generator $x \to x+\frac{\h}{q-1}$. It is easy to see that the
derivative $\pa_{q,\h}$ acts on a polynomial $f(x)$ as follows
$$
\pa_{q,\h}(f(x))= \frac{f(qx+\h)-f(x)}{(q-1)x+\h}.
$$

Also, note that besides the above analogs of the usual derivative, there are their slight modifications
$$
\tilde{\pa}_q(f(x))=\frac{f(qx)-f(\qq x)}{(q-\qq) x},\qquad \tilde{\pa}_\h(f(x))=\frac{f(x+\h)-f(x-\h)}{2\h}.
$$
These operators do not give rise to any Weyl algebra on the function space in one variable. Nevertheless, the operator $\tilde{\pa}_\h$) appears in the
frameworks of  the Weyl algebra $\W(U(u(2)_\h))$ and its commutative subalgebra considered in section 3.

We complete this section by the following observation. The q-Witt algebra contains a subalgebra looking like the enveloping algebra of the Lie algebra $sl(2)$.
Namely, consider the subalgebra generated by three  elements $e_{-1}=\pa_q$, $e_0=x\pa_q$ and $e_1=x^2\pa_q$. They are subject to the following relations
\be
e_{-1} e_0-q e_0 e_{-1}= e_{-1},\quad e_{-1} e_1 -q^2 e_1 e_{-1}=(1+q) e_0,\quad e_0 e_1-q e_1 e_0=e_1.
\label{Jackson}
\ee

This quadratic-linear algebra was considered in \cite{LS} in the frameworks of the so-called Hom-Lie algebras This notion is based on the modified Leibniz rule (\ref{modifL}).
In a similar manner the notions of  Hom-associative algebras, Hom-Poisson algebras etc., were introduced. Emphasize that the enveloping algebra of a Hom-Lie
algebra is not associative one but Hom-associative one. By contrast, our approach is based  only on the permutation relations between derivative(s) and generator(s)
of a given algebra. Thus, the algebra defined by the relations (\ref{Jackson}) is a usual associative algebra. We claim that this algebra has the good deformation
property. This property will be proven in the next section for a larger family of quadratic-linear algebras.

\section{Other roles and forms of Jacobi condition for quadratic-linear algebras}

In the previous section we presented a form of the Jacobi condition, which is useful for proving or denying the PBW property of a given quadratic-linear(-constant)
algebra. Nevertheless, in the classical case (i.e. the enveloping algebra of a usual Lie algebra is considered as such an algebra) and in some other cases, mentioned
at the end of the paper, the Jacobi condition enables one to construct the adjoint representation of a given Lie algebra. Besides, the construction of the Chevalley-Eilenberg
complex associated with this algebra is mainly based on the Jacobi identity.

In this section we discuss other forms   the Jacobi condition which are useful for  generalizing the notion of the adjoint representation and constructing an analog of the
Chevalley-Eilenberg complex on certain quadratic-linear algebras.

Let  again $U$ be a three dimensional space and $\{x,y,z\}$ be its basis. Consider the quadratic-linear algebra generated by these  generators subject to the relations
$$
xy-a yx-l_1=0,\quad yz-b zy-l_2=0,\quad zx-c xz-l_3=0,
$$
where $a,b, c\in \K$ are nontrivial constant and $l_1, l_2,l_3$ are some elements of $U$. As usual, we also consider the corresponding quadratic algebra which is obtained
by setting $l_1=l_2=l_3=0$.

It would be interesting to classify  all families $(a,b,c, l_1, l_2, l_3)$ such that the Jacobi-PP condition for the corresponding quadratic-linear algebras is valid. We restrict
ourselves to two examples.

The first example is $sl(2)$ like. We assume that
$$
l_1=kx,\quad l_2=lz,\quad  l_3=my,\quad k,l,m \in \K,\,\, klm\not=0.
$$
It is easy to see that the Jacobi-PP condition is valid iff $b=a$ and $l=k$. Thus, the relations on the generators of the corresponding algebra become
\be
xy-ayx=kx,\,\, yz-azy=kz,\,\, zx-cxz=my.
\label{rre}
\ee
This algebra can be treated as a multiparameter deformation of the commutative algebra $\Sym(sl(2))$. It can be easily seen that the algebra defined by (\ref{Jackson}) is a
particular case of the associative algebra defined by (\ref{rre}). Indeed, by identifying $e_{-1}=x$, $e_0=y$, $e_1=z$ we get (\ref{Jackson}) if in (\ref{rre}) we put $a=q$,
$c=q^{-2}$, $m=-(q^{-1}+q^{-2})$.

Since the quadratic algebra corresponding to the quadratic-linear algebra defined by (\ref{rre}) is Koszul, the latter algebra mets the  PBW property. Note that by using the
triangular structure of the quadratic-linear algebra, it is possible to define analogs of the Verma modules over it.

The second example is $su(2)$ like. We assume that
$$
l_1=kz,\quad l_2=lx,\quad l_3=my,\quad k,l,m \in \K, \,\, klm\not=0.
$$
The Jacobi-PP condition is fulfilled iff $a=b=c$. Also, by a change of a basis (over the field $\K=\C$) we can get that $k=l=m=1$. Thus, we assume that the generators are
bound by the following relations
\be
xy-ayx=z,\quad yz-azy=x,\quad zx-axz=y.
\label{sy}
\ee
Since the corresponding  quadratic algebra  is Koszul, we conclude that the algebra defined by (\ref{sy}) has the good deformation property. Thus, by introducing a parameter
$\h$ in front of the right hand side, we get a two parameter deformation of the algebra $\Sym(V)$ where $V$ is three-dimensional vector space.

Note that similar quadratic-linear-constant algebras (but with different right hand side of the relations) appear in quantization of the Poisson structures on Painlev\'e monodromy
manifolds (see \cite{MR}).

Now, we pass to considering the other forms of the Jacobi condition. Consider a quadratic algebra $\A=T(U)/\langle I \rangle$, $I\subset U^{\ot 2}$ and a quadratic-linear one
$\A_{[\,,\,]}$ (we use the notations of the previous section). Assume that the Jacobi-PP condition for the bracket $[\,,\,]:\,I \to U$ is valid. However, this condition cannot be written
in the form
\be
[\,,\,]\circ [\,,\,]_{12}=[\,,\,]\circ [\,,\,]_{23}\quad {\rm on }\quad I^{(3)}=I\ot U\bigcap U\ot I,
\label{Ja}
\ee
because the images of the operators  $[\,,\,]_{12}$ and  $[\,,\,]_{23}$ acting on $I^{(3)}$do not belong in general to $I$ (but their difference does by assumption). Thus, the
sides of (\ref{Ja}) are not well defined separately.

In order to make this object more similar to a usual Lie algebra  we assume that in the space $U^{\ot 2}$ there is a complementary subspace $I_+$ (playing the role of the
symmetric subspace) where the bracket acts trivially. Consequently, the bracket becomes well-defined on the whole space $U^{\ot 2}$. Thus, we have the following data
$(U, I, I_+, [\,,\,]:\,U^{\ot 2}\to U)$ with complementary subspaces $I\oplus I_+=U^{\ot 2}$ and such that the image of the subspace $I_+$ under the map $[\,,\,]$ is trivial.
Thus, the both sides of (\ref{Ja}) are well defined.

\begin{definition} \rm
We say that the data $(U, I, I_+, [\,,\,]:\,U^{\ot 2}\to U)$ satisfy the strong Jacobi condition if the Jacobi-PP condition is valid for the corresponding quadratic-linear algebra,
and the both sides of (\ref{Ja}) are trivial.
\end{definition}

Note that the strong Jacobi identity enables us to define an analog of the Chevalley-Eilenberg complex composed of the terms
$$
I^{(k)}=I\ot U^{\ot(k-2)} \bigcap U\ot I\ot U^{\ot(k-3)} \bigcap\dots \bigcap U^{\ot(k-3)} \ot I\ot U \bigcap U^{\ot(k-2)} \ot I
$$
with the differential $d=[\,,\,]_{12}$. The relation $d^2=0$ follows immediately from the fact that $[\,,\,] [\,,\,]_{12}=0$. Observe that elements of the subspaces
$I^{(k)}\subset U^{\ot k}$ are analogs of the space of totaly skew-symmetric elements. For this reason, even in the classical case, it suffices to apply the bracket only to two
first terms instead of employing the usual formula; the results of applying this operator $d$ and the usual Chevalley-Eilenberg operator differ by a nontrivial factor.

Now, we go back to the above examples and examine the problem whether a given quadratic-linear algebra  can be completed with convenient subspaces $I_+\subset U^{\ot 2}$
such that the new data $(U, I, I_+, [\,,\,]:\,U^{\ot 2}\to U)$ meets the strong Jacobi condition.

First, consider the $sl(2)$ like algebra. By a straightforward computation it is possible to check that the strong Jacobi condition is met if $I_+$ contains the elements $y^2$
and $cxz+azx$.

As for the $su(2)$ like example, we have that  the strong Jacobi identity is valid for the extended data iff the subspace $I_+$ contains the term $x^2+ y^2+ z^2$. In particular,
we can put
\be
I_+=\span(x^2,\, y^2,\, z^2,\,xy+\al yx,\, yz+\al zy,\, zx+\al xz),\quad \al\in \K, \,\,\al\not=0.
\label{plus}
\ee

Let us discuss now a form of the Jacobi condition for quadratic-linear algebras enabling one to construct an analog of the adjoint representation.

\begin{definition} \rm
We say that the data $(U, I, I_+, [\,,\,]:U^{\ot 2}\to U)$ with complementary $I$ and $I_+$, such that $[\,,\,] I_+=0$, is an almost Lie algebra,  if there exists a
nontrivial number $p\in \K$ such that the map $p\,[\,,\,]$ defines left and right representations of the quadratic-linear algebra defined by
\be
T(U)/\langle I-[\,,\,]I\rangle \label{ql}.
\ee
\end{definition}
Emphasize that in the classical case the normalizing factor $p$ equals 2 since the quadratic-linear algebra (\ref{ql}) differs from the usual enveloping algebra (see remark at the
end of the section). In general, it has to be found.

Let us turn to the above examples and examine the following problem: for what values of parameters  entering the defining relations of these algebras they acquire
almost Lie algebra structures.

First, consider the $su(2)$ like example. Let us assume that the subspace $I_+$ is given by (\ref{plus}). Then we get a bracket defined on the whole space $U^{\ot 2}$ and
having the following multiplication table
$$
[x,x]=0,\quad [x,y]=\frac{\al z}{\ga},\quad [y,x]=\frac{-z}{\ga},\quad {\rm c.p.}
$$
where $\ga=a+\al$.

Now, we are able to define the left action (denoted $\tri$) of the space $U$ onto itself which is multiple (with the factor $p$) of the above bracket action. We have
$$
x\tri x =0,\quad x\tri y=\frac{p\al z}{\ga},\quad x\tri z=\frac{-py}{\ga}
$$
and so on.

Thus, we can represent $x, y$ and $z$ as operators acting in the space $U=\span(x,y,z)$. We have to check that the defining relations of the quadratic-linear algebra
in question are preserved by this representation. This implies the following  relations on the parameters $a$, $\al$ and $p$:
$$
p=\al\, \gamma,\quad a \al^2 p=\gamma,\quad \gamma= a+\al.
$$
Treating $\al$ as an independent parameter we can express other parameters in terms of $\al$ as follows
$$
a=\al^{-3},\quad  \ga=\al+\al^{-3},\quad p=\al^2+\al^{-2}.
$$
It is somewhat straightforward checking that the same parameters are convenient for defining the right adjoint action of the quadratic-linear algebra in question.
In conclusion, we get a family of almost Lie algebra structures, parameterized by $\al$.

Note that the family of the data satisfying the strong Jacobi condition  is larger since the parameters $a$ entering (\ref{sy}) and $\al$ entering (\ref{plus}) are not related.

Now, pass to the $sl(2)$ type example. We assume that the subspace $I_+$ has the following form
$$
I_+=\span(x^2, \,y^2, \,z^2, \,xy+\al yx, \,yz+\beta zy, \, zx+\beta xy).
$$
By a change of the basis we can get that $k=1, l=2$ (see (\ref{rre})). By tedious but straightforward computations it can be shown that the data
$(U, I, I_+, [\,,\,]:U^{\ot 2}\to U)$ determines an almost Lie structure  iff $a=b=\al=\beta=1,\, p=2$, i.e. when it corresponds to  the usual $sl(2)$ Lie structure.

Other examples interesting from the viewpoint of the different forms of the Jacobi condition arise from braidings. In the early 80's one of the authors (D.G.) introduced the
notion of a generalized Lie algebra associated with involutive symmetries (see \cite{G1, G2}). A $gl(m)$ type example can be constructed as follows. Let  $R:\vv\to \vv$ be
a skew-invertible (see section 4) involutive symmetry. Then it can be extended up to an involutive symmetry 
$$\Ren:\End(V)^{\ot 2}\to \End(V)^{\ot 2}.$$ 
Besides, in the space $End(V)$ there is the usual product (composition) of endomorphisms
$$
\End(V) \ni X, Y  \mapsto X\circ Y\in \End(V).
$$

Introduce the following bracket
$$
[\,,\,]: \,\End(V)^{\ot 2}\to \End(V),\quad [X,Y]=X\circ Y- \circ \Ren(X,  Y).
$$

Consider the quadratic-linear algebra defined by the relations
\be
X_i \, X_j- \Ren(X_i,  X_j)= [X_i, X_j]
\label{enve}
\ee
where $\{X_i\}$ is a basis of the space $\End(V)$. We claim that for this filtered algebra the Jacobi-PP condition is valid. Moreover, the corresponding quadratic algebra
$$
T(\End(V))/\langle  I \rangle,
$$
where $I=\span(X_i \, X_j- \Ren(X_i,  X_j))$ is Koszul. This ensures the PBW property for the algebra defined by (\ref{enve}).

In contrast with the above examples where the bracket was only defined on a subspace $I\subset U^{\ot 2}$, now we have defined this bracket on the whole space
$U^{\ot 2}$, $U=\End(V)$. Consequently, we have no choice for the subspace $I_+$,  we put
$$
I_+= \span(X_i \, X_j+ \Ren(X_i, X_j)).
$$
It is not difficult to see that  the bracket above acts on this space trivially.  We claim that for the data $(U, I, I_+, [\,,\,]:U^{\ot 2}\to U)$ the strong Jacobi identity is also
valid. Moreover, this data defines an almost Lie algebra structure with $p=1$.

\begin{remark} \rm
Here $p=1$ since $[\,,\,](X_i\ot X_j -\Ren(X_i,  X_j))=2[X_i,X_j]$. Also note that for the Jacobi-PP condition, as well as for the strong Jacobi condition, the normalizing factor
$p$ does not matter.
\end{remark}

Note that the algebra defined by (\ref{enve}) can be written in the form of the modified RE algebra with an involutive $R$. So, as was noticed in section 4, this algebra has a
representation theory similar to that of the Lie algebra $gl(m)$ (or to that of $gl(m|n)$ depending on the initial symmetry $R$). Also, the spaces $I_+^{(k)}$ (resp., $I^{(k)}$)
of totaly  symmetric (resp., skew-symmetric) elements can be introduced via the projectors of symmetrization (resp., skew-symmetrization) naturally associated with the
operator $\Ren$.

Besides, considering the subspace of $\End(V)$ consisting of the traceless elements where the  {\em trace } is associated with $R$ --- the so-called braided trace, we can get
a $sl$-type data with similar properties.

If $R$ is a skew-invertible Hecke symmetry, a similar construction mutatis mutandis can be  also defined. But the role of the operator $\Ren$ is played by another operator
$Q:\End(V)^{\ot 2}\to \End(V)^{\ot 2}$ which coincides with $\Ren$ as $q\to 1$ (i.e. as the symmetry becomes involutive). It equals to the operator $Q$ from section 4 (up
to a change of $R$ by $R^{-1}$). Explicitly, the space $I$ is just the left hand side of (\ref{mRE}) whereas $I_+$ is defined by (\ref{syme}). As for the corresponding filtered
quadratic-linear algebra, it is defined by (\ref{mRE}). As usual, we can define a bracket by assuming that it kills $I_+$ and  maps the left hand side of (\ref{mRE}) to its right
hand side. Alternatively the bracket can be introduced via the product $\circ$. These two approaches lead to brackets which differ by a factor.

The point is that for this quadratic-linear algebra the Jacobi-PP condition is valid (see \cite{G3}). Since the quadratic algebra $T(\End(V))/\langle I \rangle$ is Koszul (at least
if $R$ is a deformation of an involutive symmetry and $q$ is generic) the modified RE  algebra meets the PBW property. Also, the corresponding data
$(U, I, I_+, [\,,\,]:U^{\ot 2}\to U)$, $U=\End(V)$, defines an almost Lie algebra structure. However, we do not know whether the strong Jacobi condition is valid for it.

Let us sum up. In section 3 we introduced the notion of the Weyl algebra related to an enveloping algebra of a Lie algebra. In section 4 we generalized this notion to the case
when the role of the enveloping algebra is played by a modified RE algebra. Whether it is possible to generalize this notion to the case when the enveloping algebra is replaced
by a quadratic-linear algebra?  Here we suggest a version of such a generalization.

Let $U$ be a finite-dimensional vector space with a basis $\{x_i\}$. As usual, consider a quadratic algebra $\A=T(U)/\langle I \rangle$, where  $I$ is a  subspace of $U^{\ot 2}$,
and introduce its quadratic-linear deformation $\A_{[\,,\,]}=T(U)/\langle I-[\,,\,] I \rangle$, where $[\,,\,]:\,I\to U$ is a linear map (called bracket). Also, consider  the quadratic
algebra $\B = T(U^*)/\langle J \rangle$ where $J\subset (U^*)^{\ot 2}$. Let $\pa^i$ be the basis of the space $U^*$, dual to the basis $\{x_i\}$:  $\langle \pa^i, x_j\rangle=\de_i^j$.
Let  some permutation relations of the form
\be
\pa^i x_j- \al_{j,l}^{i, k} x_k \pa^l=b_{j,k}^i \pa^k +\de_j^i
\label{perel}
\ee
be given. The space spanned by the left hand side of these elements will be denoted $K$. Assume that these permutation relations are compatible with the algebras $\A_{[\,,\,]}$
and $\B$, i.e. modulo these permutation relations any element of the product $\B\ot \A_{[\,,\,]}$ can be converted into an element of  $\A_{[\,,\,]}\ot \B$.

\begin{definition} \rm
The algebra $\W(\A_{[\,,\,]})$ generated by the algebra $\A_{[\,,\,]}$ and the algebra $\B$ whose generators are subject to the permutation relations
(\ref{perel}), is called the Weyl algebra corresponding to the algebra $\A_{[\,,\,]}$ if
\begin{enumerate}
\item
On converting the generators $\pa^j$ (called the partial derivatives) of the algebra $\B$ into operators on $\A_{[\,,\,]}$ by the same method as above (namely, with the help of
the counit $\varepsilon$), we get a representation of the algebra $\B$.
\item
The subspace
$$
I\oplus J \oplus K \subset (U\oplus U^*)^{\ot 2}
$$
endowed with the bracket, which equals $[\,,\,]$ on $I$, trivial on $J$ and defined by (\ref{perel}) on $K$, satisfies the Jacobi condition in the form of \cite{BG}.
\end{enumerate}
\end{definition}

The form of the Jacobi condition presented in \cite{BG} is a generalization of the Jacobi-PP condition covering the case of quadratic-linear-constant algebras. As was shown in
\cite{BG}, if for a given quadratic-linear-constant algebra the corresponding quadratic algebra is Koszul, then in the former algebra the PBW property is valid.

We point out that the Weyl algebras defined in section 3 and the braided Weyl algebra corresponding to a modified RE algebra are covered by this definition. It also covers the
algebra defined by (\ref{perr}). In this case the subspaces $I$ and $J$ are trivial and $K$ is generated by the left hand side of (\ref{perr}). However, we do not know whether it
is possible to define the Weyl algebras either on the Jackson algebra (\ref{Jackson}) or on $sl(2)$ and $su(2)$ like algebras considered in this section.

Nevertheless, before constructing a Weyl algebra on a given quadratic-linear algebra it is reasonable to construct such an algebra on the corresponding quadratic algebra.
The notion of a Weyl algebra corresponding to a quadratic algebra can be obtained from the  definition above by assuming the initial bracket to be trivial. There are known
numerous attempts to define partial derivatives on quadratic algebras related to the Quantum Groups. We refer the reader to the paper \cite{GS2} where some references
are given.

\end{document}